\documentclass{article}

\title{THE PARTIAL DIFFERENTILAL COEFFICIENTS FOR 
THE SECOND WEIGHTED BARTHOLDI ZETA FUNCTION OF A GRAPH}
\author{
Shigeki MATSUTANI \\ 
College of Science and Engineering, \\  
Kanazawa University, \\ 
Kanazawa, Ishikawa 920-1192, JAPAN \\ 
s-matsutani@se.kanazawa-u.ac.jp \\ 
Hideo MITSUHASHI \\
Department of Applied Informatics, \\
Faculty of Science and Engineering, \\
Hosei University, \\
Koganei, Tokyo 184-8584, JAPAN  \\
e-mail: hmitsu@hosei.ac.jp \\ 
Hideaki MORITA \\ 
Division of System Engineering for Mathematics, \\
Muroran Institute of Technology,\\ 
Muroran, Hokkaido 050-8585, JAPAN \\
e-mail: morita@mmm.muroran-it.ac.jp \\ 
Iwao SATO  \\ 
National Institute of Technology, Oyama College, \\ 
Oyama, Tochigi 323-0806, JAPAN \\
e-mail: isato@oyama-ct.ac.jp \\ 
}
 
\begin{document}
 \maketitle

Running head title:

The partial derivative of the second weighted Bartholdi zeta function  

\vspace{5mm}

The address for manuscript correspondence:

Iwao Sato 

Oyama National College of Technology, 
Oyama, Tochigi 323-0806, JAPAN

Tel: +81-285-20-2176

Fax: +81-285-20-2880

E-mail: isato@oyama-ct.ac.jp

\clearpage

\begin{abstract}
We consider the second weighted Bartholdi zeta function of a graph $G$, and present 
weighted versions for the results of Li and Hou's on the partial derivatives of the determinant part 
in the determinant expression of the Bartholdi zeta function of $G$. 
Furthermore, we give a formula for the weighted Kirchhoff index of 
a regular covering of $G$ in terms of that of $G$.  
\end{abstract}

\vspace{5mm}

 
\vspace{5mm}

\clearpage

\section{Introduction}

\subsection{The Ihara zeta function, the complexity and the Kirchhoff index of a graph}

Graphs and digraphs treated here are finite.
Let $G$ be a connected graph and $D_G$ the symmetric digraph 
corresponding to $G$. 
Set $D(G)= \{ (u,v),(v,u) \mid uv \in E(G) \} $. 
For $e=(u,v) \in D(G)$, set $u=o(e)$ and $v=t(e)$. 
Furthermore, let $e^{-1}=(v,u)$ be the {\em inverse} of $e=(u,v)$. 

The {\em Ihara(-Selberg) zeta function} of $G$ is defined by  
\[
{\bf Z} (G,t)= \prod_{[C]} (1- t^{ \mid C \mid } )^{-1} , 
\]
where $[C]$ runs over all equivalence classes of prime, reduced cycles 
of $G$ and $|C|$ is the length of a prime, reduced cycle $C$.  
Ihara [14] defined Ihara zeta functions of graphs, and showed that the 
reciprocals of Ihara zeta functions of regular graphs are explicit polynomials. 
The Ihara zeta function of a regular graph $G$ associated with a unitary 
representation of the fundamental group of $G$ was developed 
by Sunada [23,24]. 
Hashimoto [12] generalized Ihara's result on the Ihara zeta function of 
a regular graph to an irregular graph, and showed that its reciprocal is 
again a polynomial by a determinant containing the edge matrix. 
Bass [3] presented another determinant expression for the Ihara zeta function 
of an irregular graph by using its adjacency matrix. 

Let $G$ be a connected graph with $n$ vertices $v_1, \cdots ,v_n $ and $m$ edges. 
Then the {\em adjacency matrix} ${\bf A} (G)=(a_{ij} )$ is 
the square matrix such that $a_{ij} =1$ if $v_i$ and $v_j$ are adjacent, 
and $a_{ij} =0$ otherwise. 
The {\em degree} $\deg {}_G v_i $ of a vertex $v_i $ of $G$ is the number of vertices 
adjacent to $v_i $ in $G$. 
Let ${\bf D} =( d_{ij} )$ be the diagonal matrix with 
$d_{ii} =deg {}_G \  v_i $, and ${\bf Q} = {\bf D} -{\bf I} $.

Bass [3] proved the following result for the Ihara zeta function.

\newtheorem{theorem}{Theorem}
\begin{theorem}[Bass] 
Let $G$ be a connected graph with $n$ vertices and $m$ edges. 
Then the reciprocal of the Ihara zeta function of $G$ is given by 
\[
{\bf Z} (G,t) {}^{-1} =(1- t^2 )^{m-n} \det( {\bf I}_n  -t {\bf A} (G)+ t^2 {\bf Q} ) .  
\]
\end{theorem}

The {\em complexity} $ \kappa (G)$($=$ the number of spanning trees in $G$) 
of a connected graph $G$ is closely related to the Ihara zeta function of $G$. 
The complexities for various graphs were given in [4,6,8].
Hashimoto expressed the complexity of a regular graph as a limit 
involving its zeta function in [12]. 
For an irregular graph $G$, Hashimoto [13] and Northshield [21] gave 
the value of $(1-t)^{-r} {\bf Z} (G, t )^{-1} $ at $t=1$ in term of 
the complexity of $G$, where $r$ is the Betti number of $G$.

\begin{theorem}[Hashimoto; Northshield]
For any finite graph $G$ such that $r>1$, we have
\[
(1-t)^{-r} {\bf Z} (G, t )^{-1} \mid {}_{t=1} = 2^r \chi (G) \kappa (G) ,
\]
where $ \chi (G)=1-r$ is the Euler number of $G$. 
\end{theorem}

For a connected graph $G$, let 
\[
f_G (t)= \det( {\bf I} -t {\bf A} (G)+ t^2 {\bf Q} ) . 
\]

For a connected graph $G$, Northshield [21] showed that the complexity 
of $G$ is given by the derivative of the above function.

\begin{theorem}[Northshield]
For a connected graph $G$, 
\[
f^{ \prime }_G (1)=2(m-n) \kappa (G ) , 
\]
where $n= \mid V(G) \mid $ and $m= \mid E(G) \mid $. 
\end{theorem}

Let $G$ be a connected graph with $n$ vertices $v_1 , \ldots , v_n $ and $m$ edges. 
Furthermore, let $d_i =d_{v_i} = \deg {}_G v_i (1 \leq i \leq n)$, and let ${\bf L} = {\bf D} -{\bf A} (G)$ 
be the {\em Laplacian} of $G$. 
Klein and Randic [16] defined the {\em resistance distance} between $v_i $ and $v_j $ in $G$ as follows: 
\[
r_{ij} = \frac{ \det ( {\bf L}^{(ij)} )} { \kappa (G)} , 
\]
where ${\bf L}^{(ij)} $ is the matrix obtained from ${\bf L} $ by deleting its $i$ th and $j$ th rows and columns. 
If $i=j$, then we set 
\[
r_{ij} =0 . 
\] 
Kline and Randic [16] introduced the resistance distance between two vertices as the effective resistance 
between two vertices when $G$ is regarded as an electric network with a resistor of 1 ohm placed on each edge. 
In fact, Bapat et al [1] proved that the resistance distance for simple connected graphs can be calculated using 
the Laplacian. 
 
The {\em Kirchhoff index} $Kf$ of $G$ is defined by using the resistance distances as follows(see [5,11,17]): 
\[
Kf= \sum_{1 \leq i<j \leq n} r_{ij} . 
\]
The Kirchhoff index of a graph is expressed by the spectra of its Laplacian(see [5,6,16]): 
\[
Kf=n \sum^n_{i=2} \frac{1}{ \mu {}_i } , 
\]
where $Spec( {\bf L} )= \{ \mu {}_1 =0, \mu {}_2 , \ldots , \mu {}_n \} $ is the set of all eigenvalues of ${\bf L} $. 
Furthermore, Chen and Zhang [7] defined the {\em multiplicative Kirchhoff index} $Kf^*$ of $G$ as follows: 
\[
Kf^* = \sum_{1 \leq i<j \leq n} d_i d_j r_{ij} . 
\] 
In [10], the {\em additive Kirchhoff index} $Kf^+$ of $G$ was defined as follows: 
\[
Kf^+ = \sum_{1 \leq i<j \leq n} ( d_i + d_j ) r_{ij} . 
\]

Somodi [22] introduced a new Kirchhoff index $Kf^z $ of a graph by using its resistance distances: 
\[
Kf^z =Kf^z (G)= \sum_{1 \leq i<j \leq n} (d_i -2)(d_j -2) r_{ij}  . 
\]
Note that $Kf^z = Kf^* -2Kf^+ +4Kf$. 

Somodi [22] showed that the new Kirchhoff index $Kf^z (G)$ of $G$ is given by the second derivative of $f_G (u)$.

\begin{theorem}[Somodi]
Let $G$ be a connected graph with $n$ vertices and $m$ edges. 
Suppose that the minimum degree $\delta (G)$ is not less than two. 
Then 
\[
f^{\prime \prime }_G (1)=2(Kf^z +2mn-2n^2 +n) \kappa (G) . 
\]
\end{theorem}

\subsection{The weighted complexity and the weighted Kirchhoff index of a graph}

Let $G$ be a connected graph and $V(G)= \{ v_1 , \cdots , v_n \}$. 
Then we consider an $n \times n$ matrix 
${\bf W} =( w_{ij} )_{1 \leq i,j \leq n }$ with $ij$ entry 
the complex variable $w_{ij}$ if $( v_i , v_j ) \in D(G)$, 
and $w_{ij} =0$ otherwise. 
The matrix ${\bf W} = {\bf W} (G)$ is called the 
{\em weighted matrix} of $G$.
Furthermore, let $w( v_i , v_j )= w_{ij}, \  v_i , v_j \in V(G)$ and 
$w(e)= w_{ij}, e=( v_i , v_j ) \in D(G)$.

In the case that $w$ is symmetric, i.e., $w(e^{-1} )=w(e)$ 
for each $e \in D(G)$, $w: D(G) \longrightarrow {\bf C} $ is considered as a symmetric function from $E(G)$ to ${\bf C} $. 
Then, for a spanning tree $T$ of $G$, let 
\[
w(T)= \prod_{e \in E(T)} w(e) . 
\]
Furthermore, let 
\[
\kappa {}_w (G)= \sum_{T} w(T) ,  
\]
where $T$ runs over all spanning tress of $G$. 
Then this sum is called the {\em weighted complexity} of $G$. 

Now, let 
\[
f_G (w,t)= \det ( {\bf I}_n -t {\bf W} + ( {\bf D}_w - {\bf I}_n ) t^2 ) , 
\]
where $n= \mid V(G) \mid $ and ${\bf D}_w =( d_{ij} )$ is the diagonal matrix 
with $d_{ii} = \sum_{o(e)= v_i } w(e)$, $V(G)= \{ v_1 , \cdots , v_n \} $.  

Mizuno and Sato [19] showed the following result for the weighted complexity of a graph.

\begin{theorem}[Mizuno and Sato]
\[
\kappa {}_w (G) =\frac{1}{2(w(G)-n)} f^{\prime} (w,1) , 
\]
where $w(G)= \sum_{vw \in E(G)} w(v,w)$. 
\end{theorem} 

When $w=1$, i.e., $w(e)=1$ for each $e \in D(G)$, we obtain Theorem 3. 

Let $G$ be a connected graph with $n$ vertices $v_1 , \cdots , v_n $, 
$w: D(G) \longrightarrow {\bf C} $ a symmetric weight function and ${\bf W} $ 
the weighted matrix of $G$ corresponding to $w$. 
Furthermore, let ${\bf L} ={\bf D}_w - {\bf W} $ be the weighted Laplacian of $G$. 
Set $d^w_i = \sum_{o(e)= v_i } w(e)$ for each $i=1, \ldots , n$.   
For $1 \leq p \neq q \leq n$, let 
\[
r_{pq} = \frac{\det ({\bf L}^{(pq)} )}{\kappa {}_w (G)} , 
\]
where $ {\bf L}^{(pq)} =  {\bf L} (pq)$ is the submatrix of ${\bf L} $ obtained from deleting the $p, q$ th rows and 
the $p, q$ th columns. 
Then we define three weighted Kirchhoff indices of $G$ as follows: 
\[
Kf_w = \sum_{1 \leq i<j \leq n} r_{ij} ,
\]
\[
Kf^*_w = \sum_{1 \leq i<j \leq n} d^w_i d^w_j r_{ij}  
\] 
and 
\[ 
Kf^+_w = \sum_{1 \leq i<j \leq n} ( d^w_i + d^w_j ) r_{ij} . 
\]      
Furthermore, another weighted Kirchhoff index $Kf^z_w =Kf^z_w (G)$ of $G$ is defined as follows: 
\[
Kf^z_w =Kf^z_w (G)= \sum_{1 \leq p<q \leq n} (d^w_p -2)(d^w_q -2) r_{pq} .  
\]
Note that $Kf^z_w = Kf^*_w -2Kf^+_w +4Kf_w$.

The following theorem is a generalization of Theorem 4 (see [18]). 

\begin{theorem}[Mitsuhashi, Morita and Sato]
Let $G$ be a connected graph with $n$ vertices $v_1 , \cdots , v_n $, 
$w: D(G) \longrightarrow {\bf C} $ a symmetric weight function and ${\bf W} $ 
the weighted matrix of $G$ corresponding to $w$. 
Then 
\[
f^{ \prime \prime }_G (w,1)=2 \{ Kf^z_w +2w(G)n -2n^2 +n\} \kappa {}_w (G). 
\]
\end{theorem}

\subsection{Barthlodi zeta function and the Kirchhoff index of a graph}

Let $G$ be a connected graph. 
Then the {\em Bartholdi zeta function} ${\zeta} (G, u,t)$ of 
a graph $G$ is defined as follows(see [2]): 
\[
\zeta (G,u,t)= \prod_{[C]} (1- u^{ cbc(C)} t^{ \mid C \mid } )^{-1} , 
\]
where $[C]$ runs over all equivalence classes of prime cycles 
of $G$ and $cbc(C)$ is the cyclic bump count of a cycle $C$. 
Bartholdi [2] gave a determinant expression of the Bartholdi zeta 
function of a graph.

\begin{theorem}[Bartholdi] 
Let $G$ be a connected graph with $n$ vertices and $m$ unoriented edges. 
Then the reciprocal of the Bartholdi zeta function of $G$ is given by 
\[
\zeta (G,u,t )^{-1} =(1-(1-u )^2 t^2 )^{m-n} 
\det ( {\bf I}_n -t {\bf A} (G)+(1-u)( {\bf D} -(1-u) {\bf I}_n ) t^2 ) .    
\] 
\end{theorem}

Let $G$ be a connected graph with $n$ vertices $v_1 , \cdots , v_n $. 
Then Li and Hou [17] introduced the following function: 
\[
F(u,t)= \det( {\bf I}_n -t {\bf A} (G)+(1-u) t^2 ( {\bf D} -(1-u) {\bf I}_n )) . 
\]

Li and Hou [17] gave a generalization of Northshield Theorem [21].

\begin{theorem}[Li and Hou] 
Let $G$ be a connected graph with $n$ vertices and $m$ edges.   
Then 
\[ 
\kappa (G)= \frac{1}{2t^{n-2} (mt-n)} \frac{ \partial F(u,t)}{ \partial t} |_{(1-1/t,t)} , 
\] 
\[ 
\kappa (G)=- \frac{1}{2t^{n} (mt-n)} \frac{ \partial F(u,t)}{ \partial u} |_{(1-1/t,t)} . 
\]   
\end{theorem}

For $t=1$, the following result holds(see [15]).

\newtheorem{corollary}{Corollary} 
\begin{corollary}[Kim, Kwon and Lee] 
\[
\frac{ \partial F(u,t)}{ \partial t} |_{(0,1)} =- \frac{ \partial F(u,t)}{ \partial u} |_{(0,1)}
=2(m-n) \kappa (G). 
\]
\end{corollary}

Furthermore, Li and Hou [17] gave a result for the second differential coefficients of $F(u,t)$. 
Set $Kf^z (t)= Kf^* t^2 -2Kf^+ t+4Kf$.

\begin{theorem}[Li and Hou]
Let $G$ be a connected graph with $n$ vertices and $m$ edges.  
Then 
\[
\frac{ \partial {}^2 F(u,t)}{ \partial t^2} |_{(1-1/t,t)} =2t^{n-4} \kappa (G) (Kf^z (t)+nt(2mt-2n+1)) , 
\]
\[
\frac{ \partial {}^2 F(u,t)}{ \partial t \partial u} |_{(1-1/t,t)} =2t^{n-2} \kappa (G)((n-mt)(n+1)t-Kf^z (t)) ,  
\] 
\[
\frac{ \partial {}^2 F(u,t)}{ \partial u^2 } |_{(1-1/t,t)} =2 t^{n} \kappa (G) (Kf^z (t)-nt) .   
\]
\end{theorem}

In this paper, we treat the second weighted Bartholdi zeta function of a graph $G$, and present 
weighted versions for the results of Li and Hou's on the partial derivatives of the determinant part 
in the determinant expression of the Bartholdi zeta function of $G$. 

In Section 2, we consider the second weighted Bartholdi zeta function of a graph $G$, and 
express its first partial derivatives by using the weighted complexity of $G$. 
In Section 3, we express the second partial derivatives of the second weighted Bartholdi zeta function of $G$ 
by using the weighted complexity and the weighted Kirchhoff index of $G$.
As an application, we give a generalization of a theorem by Hashimoto and Northshield on 
the complexity of a graph. 
In Section 4, we give formulas for the weighted Kirchhoff index of 
a regular covering of $G$ by using the weighted Kirchhoff index of $G$. 

For a general theory of graph coverings, the reader is referred to [9].

\section{A partial differential coefficient for the second weighted Bartholdi zeta function of a graph}

Mizuno and Sato [19] defined a new zeta function of a graph by using not an infinite product but 
a determinant. 

Let $G$ be a connected graph with $n$ vertices and $m$ unoriented edges, 
and ${\bf W} = {\bf W} (G)$ a weighted matrix of $G$.
Two $2m \times 2m$ matrices ${\bf B}_w = {\bf B}_w (G)=( {\bf B}_{e,f} )_{e,f \in D(G)} $ and 
${\bf J}_0 ={\bf J}_0 (G) =( {\bf J}_{e,f} )_{e,f \in D(G)} $ 
are defined as follows: 
\[
{\bf B}_{e,f} =\left\{
\begin{array}{ll}
w(f) & \mbox{if $t(e)=o(f)$, } \\
0 & \mbox{otherwise, }
\end{array}
\right. 
\\ 
{\bf J}_{e,f} =\left\{
\begin{array}{ll}
1 & \mbox{if $f= e^{-1} $, } \\
0 & \mbox{otherwise.}
\end{array}
\right.
\] 
Then the {\em second weighted Bartholdi zeta function} of $G$ is defined by 
\[
\zeta {}_1 (G,w,u,t)= 
\det ( {\bf I}_{2m} -t ( {\bf B}_w -(1-u) {\bf J}_0 ) )^{-1} . 
\]
If $w(e)=1$ for any $e \in D(G)$, then the second weighted Bartholdi zeta function of $G$ is the 
Bartholdi zeta function of $G$(see [2]). 

The determinant expression for the second weighted Bartholdi zeta function of a graph was given by Mizuno and Sato(see [20]):

\begin{theorem}[Mizuno and Sato] 
Let $G$ be a connected graph, and 
let ${\bf W} = {\bf W} (G)$ be a weighted matrix of $G$. 
Then the reciprocal of the second weighted Bartholdi zeta function of $G$ is given by 
\[
\zeta {}_1 (G,w,u,t )^{-1} =(1-(1-u )^2 t^2 )^{m-n} 
\det ({\bf I}_n -t {\bf W} (G)+(1-u) t^2 ( {\bf D}_w -(1-u) {\bf I}_n )) , 
\]
where $n= \mid V(G) \mid $, $m= \mid E(G) \mid $. 
\end{theorem}

Let $G$ be a connected graph with $n$ vertices $v_1 , \cdots , v_n $, 
and $w : D(G) \longrightarrow {\bf C} $ be a symmetric weight function. 
Then we introduce the following function: 
\[
f_w (u,t)= f_w (G,u,t)= \det( {\bf I}_n -t {\bf W} (G)+(1-u) t^2 ( {\bf D}_w -(1-u) {\bf I}_n )) . 
\]

The following theorem is a generalization of Theorem 8.

\begin{theorem} 
Let $G$ be a connected graph with $n$ vertices $v_1 , \cdots , v_n $, 
$w: D(G) \longrightarrow {\bf C} $ a symmetric weight function and ${\bf W} $ 
the weighted matrix of $G$ corresponding to $w$. 
Then 
\[
\frac{ \partial f_w (u,t)}{ \partial t} |_{(1-1/t,t)} =2t^{n-2} (w(G)t-n)  \kappa {}_w (G) , 
\]
\[
\frac{ \partial f_w (u,t)}{ \partial u} |_{(1-1/t,t)} =-2 t^{n} (w(G)t-n)  \kappa {}_w (G) .  
\] 
\end{theorem}

{\bf Proof}.  The argument is an analogue of Somodi's method [22].  
At first,   
\[
f_w (u,t)= \det ( {\bf I}_n -t {\bf W} +(1-u)t^2 ( {\bf D}_w -(1-u) {\bf I}_n )) 
= \det ((1-(1-u)^2 t^2 ) {\bf I}_n +t {\bf L} +t((1-u)t-1) {\bf D}_w ) , 
\]
where ${\bf L} = {\bf D}_w - {\bf W} $.

Let $V(G)= \{ v_1, \cdots ,v_n \} $, and 
\[
{\bf M}_w (u,t)=(m_{ij} )=(1-(1-u)^2 t^2 ) {\bf I}_n +t {\bf L} +t((1-u)t-1) {\bf D}_w .  
\]
Furthermore, let ${\bf M} {}^t_k =( m^t_{k; ij} )$ denote the matrix ${\bf M}_w (u,t)$ 
with each entry of the $k$ th row replaced by its corresponding partial derivative 
with respect to $t$. 
Then 
\[
\frac{ \partial f_w (u,t)}{ \partial t} = \sum^n_{k=1} \det ({\bf M} {}^t_k ) . 
\] 
Here, the $(i,j)$ entry of ${\bf M} {}^t_k $ is 
\[
\begin{array}{rcl} 
\  &  & m^t_{k;ij} =(1-d^w_i t-(1-u)^2 t^2 +(1-u) t^2 d^w_i ) \delta {}_{ij} +t l_{ij} \\ 
\  &   &                \\ 
\ & + & ((2(1-u)t d^w_i -2(1-u)^2 t- d^w_i -(1-d^w_i t-(1-u)^2 t^2 +(1-u) t^2 d^w_i )) \delta {}_{ij} 
+(1-t) l_{ij} ) \delta {}_{ik}  ,  
\end{array}
\] 
where $\delta {}_{ij}$ is the Kronecker delta and ${\bf L} =( l_{ij} )$.

Thus, the $k$ th row of ${\bf M}^t_k |_{(1-1/t,t)} $ is 
\[
(0, \ldots , 0, d^w_k - \frac{2}{t} , 0, \ldots , 0)+ {\bf L} {}_k ,  
\]
where ${\bf L}_i $ is the $i$ th row of ${\bf L} $. 
Furthermore, the $i(i \neq k)$ th row of ${\bf M}^t_k |_{(1-1/t,t)} $ is ${\bf L}_i $. 
Therefore,  
\[
\det ( {\bf M} {}^t_{k} )|_{(1-1/t,t)} =( d^w_k - \frac{2}{t} ) t^{n-1} \det ({\bf L}^{(k)} )+ t^{n-1} \det ( {\bf L} ) , 
\]
where ${\bf L}^{(k)} = {\bf L} (k)$ is the submatrix from ${\bf L}^{(k)} $ obtained by deleting the $k$ th row and the $k$ th column. 
Moreover, by the Matrix-Tree Theorem, we have 
\[
\det ({\bf L}^{(k)} )= \kappa {}_w (G) . 
\]
Furthermore, we have 
\[
\det ( {\bf L} )=0 . 
\]
Therefore, it follows that 
\[
\det ( {\bf M} {}^t_{k} )|_{(1-1/t,t)} =( d^w_k - \frac{2}{t} ) t^{n-1}  \kappa {}_w (G) . 
\]
Hence, 
\[
\begin{array}{rcl}  
\  &  & \frac{ \partial f_w (u,t)}{ \partial t} |_{(1-1/t,t)} = \sum^n_{k=1} \det ( {\bf M} {}^t_{k} )|_{(1-1/t,t)} \\
\  &   &                \\ 
\ & = & \sum^n_{k=1} (t d^w_k -2) t^{n-2} \kappa {}_w (G)=2(w(G)t-n) t^{n-2} \kappa {}_w (G) .   
\end{array}
\]

Next, let ${\bf M} {}^u_k =( m^u_{k; ij} )$ denote the matrix ${\bf M}_w (u,t)$ 
with each entry of the $k$ th row replaced by its corresponding partial derivative 
with respect to $u$. 
Then 
\[
\frac{ \partial f_w (u,t)}{ \partial u} = \sum^n_{k=1} \det ({\bf M} {}^u_k ) . 
\] 
Here, the $(i,j)$ entry of ${\bf M} {}^u_k $ is 
\[
\begin{array}{rcl}  
m^u_{k;ij} & = & (1-d^w_i t-(1-u)^2 t^2 +(1-u) t^2 d^w_i ) \delta {}_{ij} +t l_{ij} \\ 
\  &   &                \\ 
\ & + & ((2(1-u) t^2 -t^2 d^w_i -(1-d^w_i t-(1-u)^2 t^2 +(1-u) t^2 d^w_i )) \delta {}_{ij} 
-t l_{ij} ) \delta {}_{ik} .  
\end{array}
\]  
Thus, the $k$ th row of ${\bf M}^u_k |_{(1-1/t,t)} $ is 
\[
(0, \ldots , 0, 2t- d^w_k t^2, 0, \ldots , 0) . 
\]
Furthermore, the $i(i \neq k)$ th row of ${\bf M}^u_k |_{(1-1/t,t)} $ is ${\bf L}_i $. 
Therefore,  
\[
\det ( {\bf M} {}^u_{k} )|_{(1-1/t,t)} =(2t- d^w_k t^2 ) t^{n-1} \det ({\bf L}^{(k)} )
=(2t- d^w_k t^2 ) t^{n-1}  \kappa {}_w (G) . 
\]
Hence, 
\[
\frac{ \partial f_w (u,t)}{ \partial u} |_{(1-1/t,t)} = \sum^n_{k=1} \det ( {\bf M} {}^u_{k} )|_{(1-1/t,t)} 
= \sum^n_{k=1} (2t- d^w_k t^2 ) t^{n-1} \kappa {}_w (G)=-2(w(G)t-n) t^{n} \kappa {}_w (G) .  
\]  
Q.E.D.

For $w={\bf 1}$, Theorem 11 implies the result of Li and Hou(Theorem 8). 
Furthermore, if $t=1$, then Theorem 11 implies Theorem 5.

\section{A second partial differential coefficient for the second weighted Bartholdi zeta function of a graph}

Let $G$ be a connected graph with $n$ vertices $v_1 , \cdots , v_n $, 
and $w : D(G) \longrightarrow {\bf C} $ be a symmetric weight function. 
Then the {\em weighted Kirchhoff index function} $Kf^z_w (t)=Kf^z_w (G,t)$ of $G$ is defined as follows: 
\[
Kf^z_w (t)=Kf^z_w (G,t)= Kf^*_w t^2 -2Kf^+_w t +4Kf_w = \sum_{1 \leq p<q \leq n} (d^w_p t-2)(d^w_q t-2) r_{pq} .  
\]

The following theorem is a generalization of Theorem 9.

\begin{theorem}
Let $G$ be a connected graph with $n$ vertices $v_1 , \cdots , v_n $, 
$w: D(G) \longrightarrow {\bf C} $ a symmetric weight function and ${\bf W} $ 
the weighted matrix of $G$ corresponding to $w$. 
Then 
\[
\frac{ \partial {}^2 f_w (u,t)}{ \partial t^2} |_{(1-1/t,t)} =2t^{n-4} \kappa {}_w (G) (Kf^z_w (t)+nt(2w(G)t-2n+1)) , 
\]
\[
\frac{ \partial {}^2 f_w (u,t)}{ \partial t \partial u} |_{(1-1/t,t)} =2t^{n-2} \kappa {}_w (G)((n-w(G)t)(n+1)t-Kf^z_w (t)) ,  
\] 
\[
\frac{ \partial {}^2 f_w (u,t)}{ \partial u^2 } |_{(1-1/t,t)} =2 t^{n} \kappa {}_w (G) (Kf^z_w (t)-nt) .   
\] 
\end{theorem}

{\bf Proof}.  The argument is an analogue of Somodi's method [22].  

At first, let $1 \leq k \neq r \leq n$. 
Then, let ${\bf M} {}^{tt}_{kk} $ denote the matrix ${\bf M}_w (u,t)$ 
with each entry of the $k$ th row replaced by its corresponding second partial derivative 
with respect to $t$. 
Furthermore, let ${\bf M} {}^{tt}_{kr} $ denote the matrix ${\bf M}_w (u,t)$ 
with each entry of the $k,r$ th rows replaced by their corresponding partial derivatives  
with respect to $t$. 
Then 
\[
\frac{ \partial {}^2 f_w (u,t)}{ \partial t^2 } 
= \sum^n_{k=1} \det ({\bf M} {}^{tt}_{kk} )+2 \sum_{1 \leq k<r \leq n} \det ({\bf M} {}^{tt}_{kr} ) . 
\]

Since the non-diagonal entries of ${\bf M} {}^{tt}_{kk} $ are a linear expression of $t$, 
all non-diagonal entries in the $k$ th row of ${\bf M} {}^{tt}_{kk} $ are 0. 
Thus, the diagonal entry in the $k$ th row of ${\bf M} {}^{tt}_{kk} $ is 
\[
2(1-u) d^w_k -2(1-u)^2 ,  
\]
and so, 
\[
\det ( {\bf M} {}^{tt}_{kk} ) |_{(1-1/t,t)} = t^{n-1} ( \frac{2}{t} d^w_k  - \frac{2}{t^2 } ) \det ( {\bf L}^{(k)} ) 
=2 t^{n-3} (t d^w_k -1) \kappa {}_w (G) . 
\]
Therefore,  
\[
\sum^n_{k=1} \det ( {\bf M} {}^{tt}_{kk} )|_{(1-1/t,t)} = \sum^n_{k=1} 2 t^{n-3} (t d^w_k -1) \kappa {}_w (G) 
=2 t^{n-3} (2w(G)t-n) \kappa {}_w (G) . 
\]

Next, let $1 \leq k \neq r \leq n$. 
Then the $k$ th row of ${\bf M}^{tt}_{kr} |_{(1-1/t,t)} $ is 
\[
(0, \ldots , 0, d^w_k - \frac{2}{t} , 0, \ldots , 0) + {\bf L}_k .  
\]
Furthermore, the $r$ th row of ${\bf M}^{tt}_{kr} |_{(1-1/t,t)} $ is 
\[
(0, \ldots , 0, d^w_r - \frac{2}{t} , 0, \ldots , 0) + {\bf L}_r .  
\] 
Thus, we have 
\[
{\bf M}^{tt}_{kr} |_{(1-1/t,t)} = {\bf L}^{[kr]} +(d^w_k - \frac{2}{t} ) {\bf S}_k +(d^w_r - \frac{2}{t} ) {\bf S}_r , 
\]
where ${\bf L}^{[kr]} $ is the matrix obtained by $t$ times all rows of ${\bf L} $ except the $k, r$ th rows, 
and 
\[
{\bf S}_k=( s_{ij} ), \ s_{ij} = \delta {}_{ij} \delta {}_{ik} . 
\]
Moreover, 
\[
\det ( {\bf L}^{[kr]} +(d^w_k - \frac{2}{t} ) {\bf S}_k +(d^w_r - \frac{2}{t} ) {\bf S}_r ) 
\]
\[
= \det ( {\bf L}^{[kr]} +(d^w_r - \frac{2}{t} ) {\bf S}_r )
+(d^w_k  - \frac{2}{t} ) \det ( {\bf L}^{[kr]} (k)+(d^w_r- \frac{2}{t} ) {\bf S}_r (k) ) . 
\]

Now, we have 
\[
\begin{array}{rcl}
\ &  & \det ( {\bf L}^{[kr]} +(d^w_r - \frac{2}{t} ) {\bf S}_r ) \\ 
\  &   &                \\ 
\ & = & \det ( {\bf L}^{[kr]} )+(d^w_r - \frac{2}{t} ) \det ({\bf L}^{[kr]} (r) )) \\
\  &   &                \\ 
\ & = & t^{n-2} \det ( {\bf L} )+(d^w_r - \frac{2}{t} ) t^{n-2} \det ({\bf L}^{(r)} ) \\
\  &   &                \\ 
\ & = & 0+(d^w_r - \frac{2}{t} ) t^{n-2} \kappa {}_w (G)=(d^w_r - \frac{2}{t} ) t^{n-2} \kappa {}_w (G). 
\end{array}
\] 
Furthermore, we have 
\[
\begin{array}{rcl}
\ &  & \det ( {\bf L}^{([kr]} (k)+(d^w_r - \frac{2}{t} ) {\bf S}_r (k)) \\ 
\  &   &                \\ 
\ & = & \det ( {\bf L}^{[kr]} (k))+(d^w_r - \frac{2}{t} ) \det ({\bf L}^{[kr]} (kr) ) \\
\  &   &                \\ 
\ & = & t^{n-2} \det ( {\bf L}^{(k)} )+(d^w_r - \frac{2}{t} ) t^{n-2} \det ({\bf L}^{(kr)} ) \\
\  &   &                \\ 
\ & = & t^{n-2} \kappa {}_w (G)+(d^w_r - \frac{2}{t} ) t^{n-2} \kappa {}_w (G) r_{kr} . 
\end{array}
\] 
Therefore, it follows that 
\[ 
\begin{array}{rcl}
\ &  & \det ({\bf M}^{tt}_{kr} ) |_{(1-1/t,t)} \\
\  &   &                \\ 
\ & = & (d^w_r - \frac{2}{t} ) t^{n-2} \kappa {}_w (G)+(d^w_k  - \frac{2}{t} )
( t^{n-2} \kappa {}_w (G)+(d^w_r - \frac{2}{t} ) t^{n-2} \kappa {}_w (G) r_{kr} ) \\ 
\  &   &                \\ 
\ & = &  t^{n-2} (d^w_k +d^w_r - \frac{4}{t} ) \kappa {}_w (G) 
+ t^{n-2} (d^w_k  - \frac{2}{t} )(d^w_r - \frac{2}{t} ) \kappa {}_w (G) r_{kr} . 
\end{array}
\] 
Hence, 
\[
\begin{array}{rcl} 
\sum_{1 \leq k<r \leq n} \det ({\bf M}^{tt}_{kr} ) |_{(1-1/t,t)} & = & 
\sum_{1 \leq k<r \leq n}  t^{n-2} (d^w_k +d^w_r - \frac{4}{t} ) \kappa {}_w (G) \\
\  &   &                \\ 
\ & + & \sum_{1 \leq k<r \leq n} t^{n-2} (d^w_k  - \frac{2}{t} )(d^w_r - \frac{2}{t} ) \kappa {}_w (G) r_{kr} . 
\end{array}
\]

But, 
\[
\sum_{1 \leq k<r \leq n} (d^w_k +d^w_r - \frac{4}{t} )=2(w(G)- \frac{n}{t} )(n-1) . 
\]
Furthermore, 
\[
\begin{array}{rcl} 
\  &  & \sum_{1 \leq k<r \leq n} (d^w_k  - \frac{2}{t} )(d^w_r - \frac{2}{t} ) r_{kr} \\ 
\  &   &                \\ 
\ & = & \sum_{1 \leq k<r \leq n} d^w_k d^w_r r_{kr} - \frac{2}{t} \sum_{1 \leq k<r \leq n} (d^w_k + d^w_r ) r_{kr}  
+ \frac{4}{t^2 } \sum_{1 \leq k<r \leq n} r_{kr} \\ 
\  &   &                \\ 
\ & = & Kf^*_w - \frac{2}{t} Kf^+_w + \frac{4}{t^2 } Kf_w \\
\  &   &                \\ 
\ & = & \frac{1}{t^2 } ( t^2 Kf^*_w -2t Kf^+_w +4 Kf_w )= \frac{1}{t^2 } Kf^z_w (t) . 
\end{array}
\] 
Therefore, it follows that 
\[
\begin{array}{rcl} 
\  &  & \frac{ \partial {}^2 f_w (u,t)}{ \partial t^2} |_{(1-1/t,t)} \\
\  &   &                \\ 
\ & = & 2 t^{n-3} (2w(G)t-n) \kappa {}_w (G) 
+2(2(w(G)- \frac{n}{t} )(n-1)+ \frac{1}{t^2 } Kf^z_w (t)) t^{n-2} \kappa {}_w (G) \\
\  &   &                \\ 
\ & = & 2t^{n-4} \kappa {}_w (G) (Kf^z_w (t)+nt(2w(G)t-2n+1)) .  
\end{array}
\]

Now, let $1 \leq i \neq j \leq n$. 
Then, let ${\bf M} {}^{tu}_{ii} $ denote the matrix ${\bf M}_w (u,t)$ with each entry of the $i$ th row replaced 
by its corresponding second partial derivative 
with respect to $t$ and $u$. 
Furthermore, let ${\bf M} {}^{tu}_{ij} $ denote the matrix ${\bf M}_w (u,t)$ 
with each entry of the $i,j$ th rows replaced by their corresponding partial derivatives  
with respect to $t$ and $u$, respectively. 
Then 
\[
\frac{ \partial {}^2 f_w (u,t)}{ \partial u \partial t} = \sum^n_{i=1} \frac{ \partial }{ \partial u} \det ({\bf M} {}^{t}_{i} ) 
= \sum^n_{i=1} \det ({\bf M} {}^{tu}_{ii} )+ \sum_{1 \leq i<j \leq n} ( \det ({\bf M} {}^{tu}_{ij} )+ \det ({\bf M} {}^{tu}_{ji} ) ) . 
\]

Since the $i$ th row of ${\bf M}^{tu}_{ii}$ is 
\[
(0, \ldots , 0, -2t d^w_i +4(1-u)t, 0,  \ldots , 0) ,   
\]
\[
\begin{array}{rcl} 
\sum^n_{i=1} \det ( {\bf M} {}^{tu}_{ii} ) |_{(1-1/t,t)} & = & \sum^n_{i=1} (-2d^w_i t+4) t^{n-1} \det ( {\bf L}^{(i)} ) \\
\  &   &                \\ 
\ & = & (-4w(G)t+4n) t^{n-1} \kappa {}_w (G) \\
\  &   &                \\ 
\ & = & 4 t^{n-1} \kappa {}_w (G) (n-w(G)t) . 
\end{array}
\]

Next, let $1 \leq i \neq j \leq n$. 
Then the $i$ th row of ${\bf M}^{tu}_{ij} |_{(1-1/t,t)} $ is 
\[
(0, \ldots , 0, d^w_i - \frac{2}{t} , 0, \ldots , 0) + {\bf L}_i .  
\]
Furthermore, the $j$ th row of ${\bf M}^{tu}_{ij} |_{(1-1/t,t)} $ is 
\[
(0, \ldots , 0, 2t- d^w_j t^2 , 0, \ldots , 0) .  
\] 
Thus, we have 
\[
\begin{array}{rcl}
\ &  & \det ({\bf M}^{tu}_{ij} )|_{(1-1/t,t)} \\
\  &   &                \\ 
\ & = & (2t- d^w_j t^2 ) \det ( {\bf L}^{[ij]} (j)+(d^w_i - \frac{2}{t} ) {\bf S}_i(j) ) \\ 
\  &   &                \\ 
\ & = & (2t- d^w_j t^2 )( \det ( {\bf L}^{[ij]} (j))+(d^w_i - \frac{2}{t} ) \det ({\bf L}^{[ij]} (ij) )) \\
\  &   &                \\ 
\ & = & (2t- d^w_j t^2 ) t^{n-2} ( \det ( {\bf L}^{(j)} )+(d^w_i - \frac{2}{t} ) \det ({\bf L}^{(ij)} )) \\
\  &   &                \\ 
\ & = & (2t- d^w_j t^2 ) t^{n-2} ( \kappa {}_w (G)+(d^w_i - \frac{2}{t} ) \kappa {}_w (G) r_{ij} ) \\
\  &   &                \\ 
\ & = & (2t- d^w_j t^2 ) t^{n-2} (1+(d^w_i - \frac{2}{t} ) r_{ij} ) \kappa {}_w (G) . 
\end{array}
\]

Similarly, we have 
\[
\det ({\bf M}^{tu}_{ji} )|_{(1-1/t,t)} =(2t- d^w_i t^2 ) t^{n-2} (1+(d^w_j - \frac{2}{t} ) r_{ij} ) \kappa {}_w (G) . 
\]
Thus, 
\[
\begin{array}{rcl}
\ &  & \sum_{1 \leq i<j \leq n} ( \det ({\bf M} {}^{tu}_{ij} )+ \det ({\bf M} {}^{tu}_{ji} ) ) \\
\  &   &                \\ 
\ & = & \sum_{1 \leq i<j \leq n} t^{n-2} \{ (2t- d^w_i t^2 )+(2t- d^w_j t^2 )-2(2- d^w_i t)(2- d^w_j t) 
r_{ij} \} \kappa {}_w (G) \\
\  &   &                \\ 
\ & = & t^{n-2} \kappa {}_w (G) \sum_{1 \leq i<j \leq n}\{ (4t-( d^w_i + d^w_j ) t^2 -2(4-2( d^w_i + d^w_j )t+ d^w_i d^w_j t^2 ) r_{ij} \} \\
\  &   &                \\ 
\ & = & t^{n-2} \kappa {}_w (G) (2n(n-1)t-2w(G)(n-1) t^2 -2(4Kf_w -2Kf^+_w t +Kf^*_w t^2 )) \\
\  &   &                \\ 
\ & = & 2 t^{n-2} \kappa {}_w (G) ((n-w(G)t)(n-1)t- Kf^z_w (t)) . 
\end{array}
\] 
Therefore, it follows that 
\[
\begin{array}{rcl}
\ &  & \frac{ \partial {}^2 f_w (u,t)}{ \partial u \partial t} |_{(1-1/t,t)} \\ 
\  &   &                \\ 
\ & = & 4 t^{n-1} \kappa {}_w (G) (n-w(G)t)+2 t^{n-2} \kappa {}_w (G) ((n-w(G)t)(n-1)t- Kf^z_w (t)) \\
\  &   &                \\ 
\ & = & 2 t^{n-2} ((n-w(G)t)(n+1)t- Kf^z_w (t)) \kappa {}_w (G) . 
\end{array}
\]

Now, we consider the second partial derivative for $f_w (u,t)$ by $u$. 
Then we have 
\[
\frac{ \partial {}^2 f_w (u,t)}{ \partial u^2 } = 
\sum^n_{i=1} \det ({\bf M} {}^{uu}_{ii} )+2 \sum_{1 \leq i<j \leq n} \det ({\bf M} {}^{uu}_{ij} ) , 
\] 
where ${\bf M} {}^{uu}_{ii} $ denote the matrix ${\bf M}_w (u,t)$ with each entry of the $i$ th row replaced 
by its corresponding second partial derivative with respect to $u$, 
and ${\bf M} {}^{uu}_{ij} $ denote the matrix ${\bf M}_w (u,t)$ with each entry of the $i,j$ th rows replaced by 
their corresponding partial derivatives with respect to $u$. 

Since the $i$ th row of ${\bf M}^{uu}_{ii}$ is 
\[
(0, \ldots , 0, -2t^2 , 0,  \ldots , 0) ,   
\]
\[
\begin{array}{rcl} 
\sum^n_{i=1} \det ( {\bf M} {}^{uu}_{ii} ) |_{(1-1/t,t)} & = & \sum^n_{i=1} -2 t^2 t^{n-1} \det ( {\bf L}^{(i)} ) \\
\  &   &                \\ 
\ & = & -2n t^{n+1} \kappa {}_w (G) . 
\end{array}
\]

Next, let $1 \leq i \neq j \leq n$. 
Then the $i$ th row of ${\bf M}^{uu}_{ij} |_{(1-1/t,t)} $ is 
\[
(0, \ldots , 0, 2t- d^w_i t^2 , 0, \ldots , 0) .  
\]
Furthermore, the $j$ th row of ${\bf M}^{uu}_{ij} |_{(1-1/t,t)} $ is 
\[
(0, \ldots , 0, 2t- d^w_j t^2 , 0, \ldots , 0) .  
\] 
Thus, we have 
\[
\begin{array}{rcl}
\ &  & \det ({\bf M}^{uu}_{ij} )|_{(1-1/t,t)} \\
\  &   &                \\ 
\ & = & (2t- d^w_i t^2 )(2t- d^w_j t^2 ) \det ( {\bf L}^{(ij)} ) t^{n-2} \\
\  &   &                \\ 
\ & = & t^n ( d^w_i t-2)( d^w_j t-2) r_{ij} \kappa {}_w(G) . 
\end{array}
\] 
Therefore, 
\[
\begin{array}{rcl}
\ &  & \sum_{1 \leq i<j \leq n} \det ({\bf M} {}^{uu}_{ij} ) |_{(1-1/t,t)} \\
\  &   &                \\ 
\ & = & \sum_{1 \leq i<j \leq n} t^n ( d^w_i t-2)( d^w_j t-2) r_{ij} \kappa {}_w(G) \\ 
\  &   &                \\ 
\ & = & t^{n} \kappa {}_w (G) \sum_{1 \leq i<j \leq n} (4-2( d^w_i + d^w_j )t+ d^w_i d^w_j t^2 ) r_{ij} \\
\  &   &                \\ 
\ & = & t^{n} Kf^z_w (t)\kappa {}_w(G) . 
\end{array}
\] 
Hence, 
\[
\frac{ \partial {}^2 f_w (u,t)}{ \partial u^2 } |_{(1-1/t,t)} =-2n t^{n+1} \kappa {}_w (G)+2 t^{n} Kf^z_w (t) \kappa {}_w(G)    
= 2 \kappa {}_w (G) t^n ( Kf^z_w (t) -nt) . 
\] 
Q.E.D.

For $w={\bf 1}$, Theorem 12 implies the result of Li and Hou(Theorem 9). 
Furthermore, if $t=1$, then Theorem 12 implies Theorem 6. 
If $w={\bf 1}$ and $t=1$, the result of Somodi(Theorem 4) is obtained from Theorem 12. 

Next, we state a weighted Kirchhoff index version of Hashimoto and Northshield Theorem.

\begin{theorem}
Let $G$ be a connected graph with $n$ vertices and $m$ edges, and 
${\bf W} (G)$ a symmetric weighted matrix of $G$. 
Then 
\[
\lim_{t \rightarrow 1} \lim_{u \rightarrow 0} \frac{(1- \sqrt{1- u^2} t)^{-m+n-1} \zeta {}_1 (G,w,u,t )^{-1} 
+2^{m-n+1} (w(G)-n) \kappa {}_w (G)}{1- \sqrt{1- u^2} t} 
\]
\[
= 2^{m-n} ( Kf^z_w +(m+n)w(G)-(m+n)n+n) \kappa {}_w (G) ,  
\]
where $ w(G)= \sum_{uv \in E(G)} w(u,v) $.
\end{theorem}

{\bf Proof}.  The argument is an analogue of Somodi's method [22]. 

At first, we have 
\[
\begin{array}{rcl}
\  &  & (1- \sqrt{1- u^2} t)^{-m+n-2} \zeta {}_1 (G,w,u,t )^{-1} \\
\  &   &                \\ 
\ & = & (1-(1- u^2 ) t^2 )^{m-n} \det ( {\bf I}_n -u {\bf W} (G)+ u^2 ({\bf D}_w - {\bf I}_n ))(1- \sqrt{1- u^2} t)^{n-m-2} \\
\  &   &                \\ 
\ & = & \frac{(1+ \sqrt{1- u^2} t)^{m-n} }{(1- \sqrt{1- u^2} t)^2 } f_w (u,t) . 
\end{array}
\] 
Note that $f_w (0,1)=0$. 
Then we have  
\[
\begin{array}{rcl}
\  &  & \lim_{t \rightarrow 1} \lim_{u \rightarrow 0} \{ (1- \sqrt{1- u^2} t)^{-m+n-2} \zeta {}_1 (G,w,u,t )^{-1} 
+ \frac{2^{m-n+1} (w(G)-n) \kappa {}_w (G) }{1- \sqrt{1- u^2} t} \} \\
\  &   &                \\ 
\ & = & \lim_{t \rightarrow 1} \lim_{u \rightarrow 0} \{ \frac{(1+ \sqrt{1- u^2} t)^{m-n} }{(1- \sqrt{1- u^2} t)^2 } f_w (u,t) 
+ \frac{2^{m-n+1} (w(G)-n) \kappa {}_w (G) }{1- \sqrt{1- u^2} t} \} \\
\  &   &                \\ 
\ & = & \lim_{t \rightarrow 1} \frac{(1+t)^{m-n} f_w (0,t)+2^{m-n+1} (w(G)-n) \kappa {}_w (G)(1-t) }{(1-t)^2 } \\
\  &   &                \\ 
\ & = & 
\lim_{t \rightarrow 1} \frac{(m-n)(1+t)^{m-n-1} f_w (0,t)+(1+t)^{m-n} \frac{ \partial f_w (u,t)}{ \partial t} |_{(0,t)} 
-2^{m-n+1} (w(G)-n) \kappa {}_w (G) }{-2(1-t)} . 
\end{array}
\]

But, substituting $t=1$, Theorem 5 implies that the numerator is equal to 
\[
0+2^{m-n+1} (w(G)-n) \kappa {}_w (G)-2^{m-n+1} (w(G)-n) \kappa {}_w (G)=0 . 
\]
Therefore, by Theorems 11 and 12, it follows that 
\[
\begin{array}{rcl}
\  &  & \lim_{t \rightarrow 1} \lim_{u \rightarrow 0} \{ (1- \sqrt{1- u^2} t)^{-m+n-2} \zeta {}_1 (G,w,u,t )^{-1} 
+ \frac{2^{m-n+1} (w(G)-n) \kappa {}_w (G) }{1- \sqrt{1- u^2} t} \} \\
\  &   &                \\ 
\ & = & \lim_{t \rightarrow 1} 1/2 \{ (m-n)(m-n-1)(1+t)^{m-n-2} f_w (0,t) \\ 
\  &   &                \\ 
\ & + & 2(m-n)(1+t)^{m-n-1} \frac{ \partial f_w (u,t)}{ \partial t} |_{(0,t)} 
+(1+t)^{m-n} \frac{ \partial {}^2 f_w (u,t)}{ \partial t^2 } |_{(0,t)} \} \\ 
\  &   &                \\ 
\ & = & 1/2 \{ (m-n)(m-n-1) 2^{m-n-2} f_w (0,1) \\ 
\  &   &                \\ 
\ & + & 2(m-n) 2^{m-n-1} \frac{ \partial f_w (u,t)}{ \partial t} |_{(0,1)} 
+ 2^{m-n} \frac{ \partial {}^2 f_w (u,t)}{ \partial t^2 } |_{(0,1)} \} \\ 
\  &   &                \\ 
\ & = &  1/2 \{ 0+2(m-n) \cdot 2^{m-n-1} \cdot 2(w(G)-n) \kappa {}_w (G)  \\ 
\  &   &                \\ 
\ & + & 2^{m-n} \cdot 2 \{ Kf^z_w (1)+2w(G)n -2n^2 +n) \} \kappa {}_w (G) \} \\. 
\  &   &                \\ 
\ & = & 2^{m-n} \{ Kf^z_w +w(G)(m+n)-(m+n)n+n \} \kappa {}_w (G) . 
\end{array}
\] 
Q.E.D.

\section{Weighted Kirchhoff indices of regular coverings}

Let $G$ be a connected graph, and 
let $N(v)= \{ w \in V(G) \mid (v,w) \in E(G) \} $ for any vertex $v$ in $G$. 
A graph $H$ is called a {\em covering} of $G$ 
with projection $ \pi : H \longrightarrow G $ if there is a surjection
$ \pi : V(H) \longrightarrow V(G)$ such that
$ \pi {\mid}_{N(v')} : N(v') \longrightarrow N(v)$ is a bijection 
for all vertices $v \in V(G)$ and $v' \in {\pi}^{-1} (v) $.
When a finite group $\Pi$ acts on a graph $G$, 
the {\em quotient graph} $G/ \Pi$ is a simple graph 
whose vertices are the $\Pi$-orbits on $V(G)$, 
with two vertices adjacent in $G/ \Pi$ if and only if some two 
of their representatives are adjacent in $G$.
A covering $ \pi : H \longrightarrow G$ is said to be
{\em regular} if there is a subgroup {\it B} of the 
automorphism group $Aut \ H$ of $H$ acting freely on $H$ such that 
the quotient graph $H/ {\it B} $ is isomorphic to $G$.

Let $G$ be a graph and $ \Gamma $ a finite group.
Then a mapping $ \alpha : D(G) \longrightarrow \Gamma $
is called an {\em ordinary voltage} {\em assignment}
if $ \alpha (v,u)= \alpha (u,v)^{-1} $ for each $(u,v) \in D(G)$.
The pair $(G, \alpha )$ is called an 
{\em ordinary voltage graph}.
The {\em derived graph} $G^{ \alpha } $ of the ordinary
voltage graph $(G, \alpha )$ is defined as follows:
\begin{quotation}
$V(G^{ \alpha } )=V(G) \times \Gamma $ and $((u,h),(v,k)) \in 
D(G^{ \alpha })$ if and only if $(u,v) \in D(G)$ and $k=h \alpha (u,v) $. 
\end{quotation}
The {\em natural projection} 
$ \pi : G^{ \alpha } \longrightarrow G$ is defined by 
$ \pi (u,h)=u, (u,h) \in V( G^{ \alpha } )$. 
The graph $G^{ \alpha }$ is called a 
{\em derived graph covering} of $G$ with voltages in 
$ \Gamma $ or a {\em $ \Gamma $-covering} of $G$.
Then, a $ \Gamma $-covering $G^{ \alpha }$ is a $ \mid \Gamma \mid $-fold
regular covering of $G$ with covering transformation group $ \Gamma $.
Furthermore, every regular covering of a graph $G$ is a  
$ \Gamma $-covering of $G$ for some group $ \Gamma $(see [9]).

In the $\Gamma $-covering $G^{ \alpha } $, set $v_g =(v,g)$ and $e_g =(e,g)$, 
where $v \in V(G), e \in D(G), g\in \Gamma $. 
For $e=(u,v) \in D(G)$, the arc $e_g$ emanates from $u_g$ and 
terminates at $v_{g \alpha (e)}$. 
Note that $ e {}^{-1}_g =( e^{-1} )_{g \alpha (e)}$. 

Let ${\bf W} = {\bf W} (G)$ be a weighted matrix of $G$. 
Then we define the {\em weighted matrix} 
$\tilde{{\bf W}} = {\bf W} ( G^{ \alpha } )
=( \tilde{w} ( u_g , v_h ))$ of $ G^{ \alpha } $ 
{\em derived from} ${\bf W}$ as follows: 
\[
\tilde{w} ( u_g , v_h ) :=\left\{
\begin{array}{ll}
w(u,v) & \mbox{if $(u,v) \in D(G)$ and $h=g \alpha (u,v)$, } \\
0 & \mbox{otherwise.}
\end{array}
\right.
\]
For $ g \in \Gamma $, let the matrix $ {\bf W} {}_g =( w^{(g)}_{uv} )$ 
be defined by
\[
w^{(g)}_{uv} :=\left\{
\begin{array}{ll}
w(u,v) & \mbox{if $ \alpha (u,v)=g$ and $(u,v) \in D(G)$, } \\
0 & \mbox{otherwise.}
\end{array}
\right.
\]

If $ {\bf M}_{1} , {\bf M}_{2} , \cdots , {\bf M}_{s} $ are 
square matrices, then let \( {\bf M}_{1} \oplus \cdots \oplus {\bf M}_{s} \) 
be the block diagonal sum of ${\bf M}_{1}, \cdots , {\bf M}_{s}$ and if 
\( {\bf M}_{1} = {\bf M}_{2} = \cdots = {\bf M}_{s} = {\bf M} \),
then we write 
\( s \circ {\bf M} = {\bf M}_{1} \oplus \cdots \oplus {\bf M}_{s} \).
The {\em Kronecker product} $ {\bf A} \bigotimes {\bf B} $
of matrices {\bf A} and {\bf B} is considered as the matrix 
{\bf A} having the element $a_{ij}$ replaced by the matrix $a_{ij} {\bf B}$.

Mizuno and Sato [20] presented a relation between $f_{G {}^{ \alpha } } ( \tilde{w} ,u)$ and $f_G (w,u)$.

\begin{theorem}[Mizuno and Sato]  
Let $G$ be a connected graph with $n$ vertices, 
$ \Gamma $ a finite group and $ \alpha : D(G) \longrightarrow \Gamma  $ 
an ordinary voltage assignment. 
Moreover, let ${\bf W} (G)$ be a symmetric weighted matrix of $G$. 
Furthermore, let $ {\rho}_{1} =1, {\rho}_{2} , \cdots , {\rho}_{k} $
be all inequivalent irreducible representations of $ \Gamma $, and 
$f_i$ the degree of $ {\rho}_{i} $ for each $i$, where 
$f_1=1$.
Then 
\[
f_{ \tilde{w} } (G {}^{ \alpha } ,u,t)= 
f_w (G,u,t) \cdot \prod^{k}_{i=2} 
\det ({\bf I}_{n f_i } -t \sum_{g \in \Gamma } 
{\rho}_{i} (g) \bigotimes {\bf W}_g +(1-u) t^2 {\bf I}_{f_i } \bigotimes ( {\bf D}_w -(1-u) {\bf I}_n ))^{f_i } .  
\] 
\end{theorem}

By Theorem 14, $f_w (G,u,t)$ divides 
$f_{ \tilde{w}} (G^{\alpha } , u,t)$.

Mizuno and Sato [19] explicitly expressed the weighted complexity of a connected 
regular covering of $G$ by using that of $G$.

\begin{theorem}[Mizuno and Sato]  
Let $G$ be a connected graph with $n$ vertices, 
$ \Gamma $ a finite group and 
$\alpha : D(G) \longrightarrow \Gamma $ an ordinary voltage assignment.
Furthermore, let ${\bf W} (G)$ be a symmetric weighted matrix of $G$. 
Furthermore, let $ {\rho}_{1} =1, {\rho}_{2} , \cdots , {\rho}_{k} $
be the irreducible representations of $ \Gamma $, and 
$f_i$ the degree of $ {\rho}_{i} $ for each $i$, where 
$f_1=1$. 
Suppose that the $ \Gamma $-covering $G^{ \alpha } $ of $G$ is 
connected. 
Then the weighted complexity of $G^{ \alpha}$ is
\[
\kappa {}_{\tilde{w}} (G^{ \alpha } )= \frac{1}{ \mid \Gamma \mid } \kappa {}_w (G) 
\prod^k_{i=2} 
\det ({\bf I}_{f_i } \bigotimes {\bf D}_w - \sum_{g \in \Gamma } { \rho }_i (g) \bigotimes {\bf W}_{g} )^{f_i} . 
\]
\end{theorem}

By Theorems 14 and 15, we explicitly express the weighted Kirchhoff index of a connected 
regular covering of $G$ by using that of $G$.

\begin{theorem}
Let $G$ be a connected graph with $n$ vertices and $m$ edges, 
$ \Gamma $ a finite group and 
$\alpha : D(G) \longrightarrow \Gamma $ an ordinary voltage assignment.
Let ${\bf W} (G)$ be a symmetric weighted matrix of $G$. 
Set $| \Gamma |=r$. 
Let $ {\rho}_{1} =1, {\rho}_{2} , \cdots , {\rho}_{k} $
be the irreducible representations of $ \Gamma $, and 
$f_i$ the degree of $ {\rho}_{i} $ for each $i$, where 
$f_1=1$. 
Suppose that the $ \Gamma $-covering $G^{ \alpha } $ of $G$ is 
connected. 
Then
\begin{enumerate} 
\item 
\[
Kf^z_{\tilde{w}} ( G {}^{ \alpha } ,t)= | \Gamma | Kf^z_w (G,t)-2| \Gamma |(w(G)t -n) 
(( | \Gamma |-1)nt- t^2 \sum^k_{i=2} f_i \frac{ \partial}{ \partial t} \log K_i (u,t) |_{(1-1/t,t)} ) . 
\] 
\item
\[
\begin{array}{rcl}
\ &  & Kf^z_{\tilde{w}} ( G {}^{ \alpha } ,t) \\
\  &   &                \\ 
\ & = & | \Gamma | Kf^z_w (G,t)+| \Gamma |(n-w(G)t)((| \Gamma |-1)nt 
- t^2 \sum^k_{i=2} f_i \frac{ \partial}{ \partial t} \log K_i (u,t) |_{(1-1/t,t)} \\ 
\  &   &                \\ 
\ & + & \sum^k_{i=2} f_i \frac{ \partial}{ \partial u} \log K_i (u,t) |_{(1-1/t,t)} ) . 
\end{array}
\]
\item 
\[
Kf^z_{\tilde{w}} ( G {}^{ \alpha } ,t)=| \Gamma | Kf^z_w (G,t)-2 | \Gamma |(w(G)t-n) 
\sum^k_{i=2} f_i \frac{ \partial}{ \partial u} \log K_i (u,t) |_{(1-1/t,t)}  . 
\] 
Here, 
\[
K_i (u,t)= \det ({\bf I}_{n f_i } -t \sum_{g \in \Gamma } 
{\rho}_{i} (g) \bigotimes {\bf W}_g +(1-u) t^2 {\bf I}_{f_i } \bigotimes ( {\bf D}_w -(1-u) {\bf I}_n )) \ (2 \leq i \leq k) 
\]
\end{enumerate}
\end{theorem}

{\bf Proof}.  By Theorem 14, we have 
\[
f_{\tilde{w} } (G {}^{ \alpha } ,u,t)= f_w (G,u,t) \cdot \prod^{k}_{i=2} 
\det ({\bf I}_{n f_i } -t \sum_{g \in \Gamma } 
{\rho}_{i} (g) \bigotimes {\bf W}_g +(1-u) t^2 {\bf I}_{f_i } \bigotimes ( {\bf D}_w -(1-u) {\bf I}_n ))^{f_i } , 
\]

Now, let 
\[
K_i (u,t)= \det ({\bf I}_{n f_i } -t \sum_{g \in \Gamma } 
{\rho}_{i} (g) \bigotimes {\bf W}_g +(1-u) t^2 {\bf I}_{f_i } \bigotimes ( {\bf D}_w -(1-u) {\bf I}_n )) \ (2 \leq i \leq k) 
\]
and 
\[
K(u,t)= \prod^{k}_{i=2} K_i (u,t)^{f_i } . 
\]
Then we have 
\[
f_{\tilde{w} } (G {}^{ \alpha } ,u,t)= f_w (G,u,t) K(u,t) . 
\]
Thus, we have 
\[
\frac{ \partial {}^2 f_{\tilde{w}} (G {}^{ \alpha } ,u,t)}{\partial t^2 } = 
\frac{ \partial {}^2 f_w (G,u,t)}{\partial t^2 } K(u,t) 
+2 \frac{ \partial f_w (G,u,t)}{\partial t} \frac{ \partial K(u,t)}{\partial t} 
+ f_w (G,u,t) \frac{ \partial {}^2 K(u,t)}{\partial t^2 } . 
\]
Since 
\[
f_w (G, 1-1/t,t)= \det (t {\bf L} )=0 , 
\]
we have 
\[
\begin{array}{rcl}
\frac{ \partial {}^2 f_{\tilde{w}} (G {}^{ \alpha } ,u,t)}{\partial t^2 } |_{(1-1/t,t)} & = &  
\frac{ \partial {}^2 f_w (G,u,t)}{\partial t^2 } |_{(1-1/t,t)} K(1-1/t,t) \\ 
\  &   &                \\ 
\ & + & 
2 \frac{ \partial f_w (G,u,t)}{\partial t} |_{(1-1/t,t)} \frac{ \partial K(u,t)}{\partial t} |_{(1-1/t,t)} .  
\end{array}
\]

But, 
\[
K(0,1)= \prod^k_{i=2} 
\det ({\bf I}_{f_i } \bigotimes {\bf D}_w - \sum_{g \in \Gamma } { \rho }_i (g) \bigotimes {\bf W}_{g} )^{f_i}  
\]
and 
\[
K(1-1/t,t)= t^{( f^2_2 + \cdots + f^2_k )n} \prod^k_{i=2} 
\det ({\bf I}_{f_i } \bigotimes {\bf D}_w - \sum_{g \in \Gamma } { \rho }_i (g) \bigotimes {\bf W}_{g} )^{f_i}  
= t^{(r-1)n} K(0,1) . 
\] 
Note that $| \Gamma |=r=1+ f^2_2 + \cdots + f^2_k $. 
By Theorems 11 and 12, we have 
\[
\begin{array}{rcl} 
\  &  & \frac{ \partial {}^2 f_{\tilde{w}} (G {}^{ \alpha } ,u,t)}{\partial t^2 } |_{(1-1/t,t)} \\ 
\  &   &                \\ 
\ & = & 2 t^{nr-4} \{ Kf^z_{\tilde{w}} ( G {}^{ \alpha } ,t)+n| \Gamma |t(2w(G {}^{ \alpha } )t-2n| \Gamma |+1) \} 
\kappa {}_{\tilde {w}} (G {}^{ \alpha } ) \\
\  &   &                \\ 
\ & = & 2 t^{n-4} \{ Kf^z_w (G,t)+nt(2w(G)t -2n +1) \} \kappa {}_w (G) K(1-1/t,t) \\ 
\  &   &                \\ 
\ & + & 4 t^{n-2} \cdot (w(G)t-n) \kappa {}_w (G)  \frac{ \partial K(u,t)}{ \partial t} |_{(1-1/t,t)} . 
\end{array}
\] 
By Theorem 15, we have 
\[
\begin{array}{rcl}
\  &  & \frac{ \partial {}^2 f_{\tilde{w}} (G {}^{ \alpha } ,u,t)}{\partial t^2 } |_{(1-1/t,t)} \\ 
\  &   &                \\ 
\ & = & \frac{2}{r} t^{nr-4} \kappa {}_{w} (G)K(0,1) \{ Kf^z_{\tilde{w}} ( G {}^{ \alpha } ,t)+2nr^2 t(w(G)t-n)+nrt \} \\
\  &   &                \\ 
\ & = & 2 t^{n-4} \{ Kf^z_w (G,t)+2nt(w(G)t -n) +nt \} \kappa {}_w (G) K(1-1/t,t) \\ 
\  &   &                \\ 
\ & + & 4 t^{n-2} \cdot (w(G)t-n) \kappa {}_w (G)  \frac{ \partial K(u,t)}{ \partial t} |_{(1-1/t,t)} . 
\end{array}
\]

But, 
\[
\begin{array}{rcl}
\ &  &\frac{ \partial K(u,t)}{\partial t} |_{(1-1/t,t)} \\
\  &   &                \\ 
\ & = & \sum^k_{i=2} \prod^k_{j=2, j \neq i} K_j (1-1/t,t)^{f_j } \cdot f_i K_i (1-1/t,t)^{f_i -1} 
\frac{ \partial K_i (u,t)}{ \partial t} |_{(1-1/t,t)} \\
\  &   &                \\ 
\ & = & \sum^k_{i=2} \prod^k_{j=2} K_j (1-1/t,t)^{f_j } \frac{f_i }{K_i (1-1/t,t)} 
\frac{ \partial K_i (u,t)}{ \partial t} |_{(1-1/t,t)} \\ 
\  &   &                \\ 
\ & = & K(1-1/t,t) \sum^k_{i=2} f_i \frac{ \partial}{ \partial t} \log K_i (u,t) |_{(1-1/t,t)} \\ 
\  &   &                \\ 
\ & = & t^{(r-1)n} K(0,1) \sum^k_{i=2} f_i \frac{ \partial}{ \partial t} \log K_i (u,t) |_{(1-1/t,t)} . 
\end{array}
\] 
Thus, 
\[
\begin{array}{rcl}
\  &  & \frac{2}{r} t^{nr-4} \kappa {}_{w} (G) K(0,1) \{ Kf^z_{\tilde{w}} ( G {}^{ \alpha } ,t)+2nr^2 t(w(G)t-n)+nrt \} \\
\  &   &                \\ 
\ & = & 2 t^{n-4} \{ Kf^z_w (G,t)+2nt(w(G)t -n) +nt \} \kappa {}_w (G) t^{(r-1)n} K(0,1) \\ 
\  &   &                \\ 
\ & + & 4 t^{n-2} \cdot (w(G)t-n) \kappa {}_w (G) t^{(r-1)n} K(0,1) 
\sum^k_{i=2} f_i \frac{ \partial}{ \partial t} \log K_i (u,t) |_{(1-1/t,t)} . 
\end{array}
\]
Therefore, it follows that 
\[
\begin{array}{rcl}
\  &  & \frac{2}{r} t^{nr-4} \{ Kf^z_{\tilde{w}} ( G {}^{ \alpha } ,t)+2nr^2 t(w(G)t-n)+nrt \} \\
\  &   &                \\ 
\ & = & 2 t^{nr-4} \{ Kf^z_w (G,t)+2nt(w(G)t -n)+nt\} \\ 
\  &   &                \\ 
\ & + & 4 t^{nr-2} \cdot (w(G)t-n) 
\sum^k_{i=2} f_i \frac{ \partial}{ \partial t} \log K_i (u,t) |_{(1-1/t,t)} . 
\end{array}
\]
Hence 
\[
Kf^z_{\tilde{w}} ( G {}^{ \alpha } ,t)= | \Gamma | Kf^z_w (G,t)-2| \Gamma |(w(G)t -n) 
(( | \Gamma |-1)nt- t^2 \sum^k_{i=2} f_i \frac{ \partial}{ \partial t} \log K_i (u,t) |_{(1-1/t,t)} ) . 
\]

Next, we have 
\[
\begin{array}{rcl}
\frac{ \partial {}^2 f_{\tilde{w}} (G {}^{ \alpha } ,u,t)}{\partial u \partial t} & = & 
\frac{ \partial {}^2 f_w (G,u,t)}{\partial u \partial t} K(u,t) 
+ \frac{ \partial f_w (G,u,t)}{\partial u} \frac{ \partial K(u,t)}{\partial t} \\ 
\  &   &                \\ 
\ & + &  \frac{ \partial f_w (G,u,t)}{\partial t} \frac{ \partial K(u,t)}{\partial u} 
+ f_w (G,u,t) \frac{ \partial {}^2 K(u,t)}{\partial u \partial t} . 
\end{array}
\] 
Then we have 
\[
\begin{array}{rcl}
\frac{ \partial {}^2 f_{\tilde{w}} (G {}^{ \alpha } ,u,t)}{\partial u \partial t} |_{(1-1/t,t)} & = & 
\frac{ \partial {}^2 f_w (G,u,t)}{\partial u \partial t} |_{(1-1/t,t)} K(1-1/t,t) 
+ \frac{ \partial f_w (G,u,t)}{\partial u} |_{(1-1/t,t)} \frac{ \partial K(u,t)}{\partial t} |_{(1-1/t,t)} \\ 
\  &   &                \\ 
\ & + &  \frac{ \partial f_w (G,u,t)}{\partial t} |_{(1-1/t,t)} \frac{ \partial K(u,t)}{\partial u} |_{(1-1/t,t)} .  
\end{array}
\]

By Theorems 11 and 12, we have 
\[
\begin{array}{rcl}
\  &  & \frac{ \partial {}^2 f_{\tilde{w}} (G {}^{ \alpha } ,u,t)}{\partial u \partial t} |_{(1-1/t,t)} \\ 
\  &   &                \\ 
\ & = & 2 t^{nr-2} \{ -Kf^z_{\tilde{w}} ( G {}^{ \alpha } ,t)+(nr-rw(G)t)(nr+1)t \} 
\kappa {}_{\tilde {w}} (G {}^{ \alpha } ) \\
\  &   &                \\ 
\ & = & 2 t^{n-2} \{ -Kf^z_w (G,t)+(n-w(G)t)(n+1)t \} \kappa {}_w (G) K(1-1/t,t) \\ 
\  &   &                \\ 
\ & - & 2 t^{n} \cdot (w(G)t-n) \kappa {}_w (G)  \frac{ \partial K(u,t)}{ \partial t} |_{(1-1/t,t)} \\
\  &   &                \\ 
\ & + & 2 t^{n-2} \cdot (w(G)t-n) \kappa {}_w (G)  \frac{ \partial K(u,t)}{ \partial u} |_{(1-1/t,t)} . 
\end{array}
\] 
By Theorem 15, we have 
\[
\begin{array}{rcl}
\  &  & \frac{ \partial {}^2 f_{\tilde{w}} (G {}^{ \alpha } ,u,t)}{\partial u \partial t} |_{(1-1/t,t)} \\ 
\  &   &                \\ 
\ & = & \frac{2}{r} t^{nr-2} \kappa {}_{w} (G)K(0,1) \{ -Kf^z_{\tilde{w}} ( G {}^{ \alpha } ,t)+(nr-rw(G)t)(nr+1)t \} \\
\  &   &                \\ 
\ & = & 2 t^{n-2} \{ -Kf^z_w (G,t)+(n-w(G)t)(n+1)t \} \kappa {}_w (G) t^{(r-1)n} K(0,1) \\ 
\  &   &                \\ 
\ & + & 2 t^{n} \cdot (n-w(G)t) \kappa {}_w (G)  \frac{ \partial K(u,t)}{ \partial t} |_{(1-1/t,t)} \\ 
\  &   &                \\ 
\ & - & 2 t^{n-2} \cdot (n-w(G)t) \kappa {}_w (G)  \frac{ \partial K(u,t)}{ \partial u} |_{(1-1/t,t)} .  
\end{array}
\]

But, 
\[
\begin{array}{rcl}
\ &  &\frac{ \partial K(u,t)}{\partial u} |_{(1-1/t,t)} \\
\  &   &                \\ 
\ & = & \sum^k_{i=2} \prod^k_{j=2, j \neq i} K_j (1-1/t,t)^{f_j } \cdot f_i K_i (1-1/t,t)^{f_i -1} 
\frac{ \partial K_i (u,t)}{ \partial u} |_{(1-1/t,t)} \\
\  &   &                \\ 
\ & = & \sum^k_{i=2} \prod^k_{j=2} K_j (1-1/t,t)^{f_j } \frac{f_i }{K_i (1-1/t,t)} 
\frac{ \partial K_i (u,t)}{ \partial u} |_{(1-1/t,t)} \\ 
\  &   &                \\ 
\ & = & K(1-1/t,t) \sum^k_{i=2} f_i \frac{ \partial}{ \partial u} \log K_i (u,t) |_{(1-1/t,t)} \\ 
\  &   &                \\ 
\ & = & t^{(r-1)n} K(0,1) \sum^k_{i=2} f_i \frac{ \partial}{ \partial u} \log K_i (u,t) |_{(1-1/t,t)} . 
\end{array}
\] 
Thus, 
\[
\begin{array}{rcl}
\  &  & \frac{2}{r} t^{nr-2} \kappa {}_{w} (G)K(0,1) \{ -Kf^z_{\tilde{w}} ( G {}^{ \alpha } ,t)+(nr-rw(G)t)(nr+1)t \} \\
\  &   &                \\ 
\ & = & 2 t^{n-2} \{ -Kf^z_w (G,t)+(n-w(G)t)(n+1)t \} \kappa {}_w (G) t^{(r-1)n} K(0,1) \\ 
\  &   &                \\ 
\ & + & 2 t^{n} \cdot (n-w(G)t) \kappa {}_w (G) t^{(r-1)n} K(0,1) 
\sum^k_{i=2} f_i \frac{ \partial}{ \partial t} \log K_i (u,t) |_{(1-1/t,t)} \\ 
\  &   &                \\ 
\ & - & 2 t^{n-2} \cdot (n-w(G)t) \kappa {}_w (G) t^{(r-1)n} K(0,1) 
\sum^k_{i=2} f_i \frac{ \partial}{ \partial u} \log K_i (u,t) |_{(1-1/t,t)} .
\end{array}
\] 
Therefore, it follows that 
\[
\begin{array}{rcl}
\  &  & \frac{2}{r} t^{nr-2} \{ -Kf^z_{\tilde{w}} ( G {}^{ \alpha } ,t)+(nr-rw(G)t)(nr+1)t \} \\
\  &   &                \\ 
\ & = & 2 t^{nr-2} \{ -Kf^z_w (G,t)+(n-w(G)t)(n+1)t \} \\ 
\  &   &                \\ 
\ & + & 2 t^{nr} \cdot (n-w(G)t) \sum^k_{i=2} f_i \frac{ \partial}{ \partial t} \log K_i (u,t) |_{(1-1/t,t)} \\ 
\  &   &                \\ 
\ & - & 2 t^{nr-2} \cdot (n-w(G)t) \sum^k_{i=2} f_i \frac{ \partial}{ \partial u} \log K_i (u,t) |_{(1-1/t,t)} .
\end{array}
\] 
Hence 
\[
\begin{array}{rcl}
\ &  & Kf^z_{\tilde{w}} ( G {}^{ \alpha } ,t) \\
\  &   &                \\ 
\ & = & | \Gamma | Kf^z_w (G,t)+| \Gamma |(n-w(G)t)((| \Gamma |-1)nt 
- t^2 \sum^k_{i=2} f_i \frac{ \partial}{ \partial t} \log K_i (u,t) |_{(1-1/t,t)} \\ 
\  &   &                \\ 
\ & + & \sum^k_{i=2} f_i \frac{ \partial}{ \partial u} \log K_i (u,t) |_{(1-1/t,t)} ) . 
\end{array}
\]

Finally, we have 
\[
\begin{array}{rcl}
\frac{ \partial {}^2 f_{\tilde{w}} (G {}^{ \alpha } ,u,t)}{\partial u^2 } & = & 
\frac{ \partial {}^2 f_w (G,u,t)}{\partial u^2 } K(u,t) 
+2 \frac{ \partial f_w (G,u,t)}{\partial u} \frac{ \partial K(u,t)}{\partial u} \\ 
\  &   &                \\ 
\ & + &  f_w (G,u,t) \frac{ \partial {}^2 K(u,t)}{\partial u^2 } . 
\end{array}
\] 
Then we have 
\[
\frac{ \partial {}^2 f_{\tilde{w}} (G {}^{ \alpha } ,u,t)}{\partial u^2 } |_{(1-1/t,t)} =  
\frac{ \partial {}^2 f_w (G,u,t)}{\partial u^2 } |_{(1-1/t,t)}  K(1-1/t,t) 
+2 \frac{ \partial f_w (G,u,t)}{\partial u} |_{(1-1/t,t)}  \frac{ \partial K(u,t)}{\partial u} |_{(1-1/t,t)}  . 
\]

By Theorems 11 and 12, we have 
\[
\begin{array}{rcl}
\  &  &\frac{ \partial {}^2 f_{\tilde{w}} (G {}^{ \alpha } ,u,t)}{\partial u^2 } |_{(1-1/t,t)} \\ 
\  &   &                \\ 
\ & = & 2 t^{nr} ( Kf^z_{\tilde{w}} ( G {}^{ \alpha } ,t)-nrt)  
\kappa {}_{\tilde{w}} (G {}^{ \alpha } ) \\
\  &   &                \\ 
\ & = & 2 t^{n} ( Kf^z_w (G,t)-nt) \kappa {}_w (G) K(1-1/t,t) \\ 
\  &   &                \\ 
\ & - & 4 t^{n} \cdot (w(G)t-n) \kappa {}_w (G) \frac{ \partial K(u,t)}{ \partial u} |_{(1-1/t,t)} . 
\end{array}
\] 
By Theorem 15, we have 
\[
\begin{array}{rcl}
\  &  &\frac{ \partial {}^2 f_{\tilde{w}} (G {}^{ \alpha } ,u,t)}{\partial u^2 } |_{(1-1/t,t)} \\ 
\  &   &                \\ 
\ & = & \frac{2}{r} t^{nr} \kappa {}_{w} (G)K(0,1) ( Kf^z_{\tilde{w}} ( G {}^{ \alpha } ,t)-nrt) \\
\  &   &                \\ 
\ & = & 2 t^{n} ( Kf^z_w (G,t)-nt) \kappa {}_w (G) t^{(r-1)n} K(0,1) \\ 
\  &   &                \\ 
\ & - & 4 t^{n} \cdot (w(G)t-n) \kappa {}_w (G) t^{(r-1)n} K(0,1) 
\sum^k_{i=2} f_i \frac{ \partial}{ \partial u} \log K_i (u,t) |_{(1-1/t,t)} .  
\end{array}
\] 
Thus, 
\[
\begin{array}{rcl}
\  &  & \frac{2}{r} t^{nr} ( Kf^z_{\tilde{w}} ( G {}^{ \alpha } ,t)-nrt) \\
\  &   &                \\ 
\ & = & 2 t^{nr} ( Kf^z_w (G,t)-nt)-4 t^{nr} \cdot (w(G)t-n) 
\sum^k_{i=2} f_i \frac{ \partial}{ \partial u} \log K_i (u,t) |_{(1-1/t,t)} .  
\end{array}
\]
Therefore, it follows that 
\[
Kf^z_{\tilde{w}} ( G {}^{ \alpha } ,t)=| \Gamma | Kf^z_w (G,t)-2 | \Gamma |(w(G)t-n) 
\sum^k_{i=2} f_i \frac{ \partial}{ \partial u} \log K_i (u,t) |_{(1-1/t,t)}  . 
\] 
Q.E.D.

In the case of $t=1$, the following result holds(see [18]).

\begin{corollary}[Mitsuhashi, Morita and Sato] 
Let $G$ be a connected graph with $n$ vertices and $m$ edges, 
$ \Gamma $ a finite group and 
$\alpha : D(G) \longrightarrow \Gamma $ an ordinary voltage assignment.
Furthermore, let ${\bf W} (G)$ be a symmetric weighted matrix of $G$. 
Set ${\bf Q}_w = {\bf D}_w - {\bf I}_n $. 
Let $ {\rho}_{1} =1, {\rho}_{2} , \cdots , {\rho}_{k} $
be the irreducible representations of $ \Gamma $, and 
$f_i$ the degree of $ {\rho}_{i} $ for each $i$, where 
$f_1=1$. 
Suppose that the $ \Gamma $-covering $G^{ \alpha } $ of $G$ is 
connected. 
Then
\[
\begin{array}{rcl}
\  &  &  Kf^z_{\tilde{w}} ( G {}^{ \alpha } )=| \Gamma | Kf^z_w (G)-2n(w(G) -n)(| \Gamma |^2 -| \Gamma |) \\ 
\  &   &                \\ 
\ & + & 2| \Gamma |(w(G)-n) \sum^k_{i=2} f_i \frac{d}{dt} 
\det ( {\bf I}_{f_i n} -t \sum^k_{g \in \Gamma } \rho {}_i (g) \bigotimes {\bf W}_g +t^2 ( f_i \circ {\bf Q}_w )) |_{t=1} . 
\end{array}
\]
\end{corollary}

{\bf Proof}.  By the first formula of Theorem 16, we have 
\[
\begin{array}{rcl}
\  &  &  Kf^z_{\tilde{w}} ( G {}^{ \alpha } )= Kf^z_{\tilde{w}} ( G {}^{ \alpha } ,1) \\
\  &   &                \\ 
\ & = & | \Gamma | Kf^z_w (G,1)+2| \Gamma |n(1-| \Gamma |)(w(G)-n) \\ 
\  &   &                \\ 
\ & + & 2| \Gamma |(w(G)-n) \sum^k_{i=2} f_i \frac{ \partial}{ \partial t} \log K_i (u,t) |_{(0,1)} .
\end{array}
\] 
But, 
\[
\begin{array}{rcl}
\  &  &  \frac{ \partial}{ \partial t} \log K_i (u,t) |_{(0,1)} = \frac{d}{dt} \log K_i (0,t) |_{t=1} \\
\  &   &                \\ 
\ & = & \frac{d}{dt} \det ( {\bf I}_{f_i n} -t \sum^k_{g \in \Gamma } \rho {}_i (g) \bigotimes {\bf W}_g 
+t^2 ( f_i \circ {\bf Q}_w )) |_{t=1} .
\end{array}
\] 
Therefore,  the result follows. 
Q.E.D.

Next, in the case of $w= {\bf 1} $, we have 
\[
w(G)=|E(G)|, \ Kf^z_{{\bf 1}} (G)=Kf^z (G), \ Kf^z_{\tilde{{\bf 1}}} ( G^{\alpha } )=Kf^z ( G^{\alpha } ), \ 
{\bf D} = {\bf D}_{{\bf 1}} 
\]
and 
\[ 
{\bf Q} = {\bf Q}_{{\bf 1}} = {\bf D} - {\bf I}_n . 
\]
For $g \in \Gamma $, let 
\[
{\bf A}_g = {\bf W}_g . 
\]

Thus,

\begin{corollary}
Let $G$ be a connected graph with $n$ vertices and $m$ edges, 
$ \Gamma $ a finite group and 
$\alpha : D(G) \longrightarrow \Gamma $ an ordinary voltage assignment.
Let $ {\rho}_{1} =1, {\rho}_{2} , \cdots , {\rho}_{k} $
be the irreducible representations of $ \Gamma $, and 
$f_i$ the degree of $ {\rho}_{i} $ for each $i$, where 
$f_1=1$. 
Suppose that the $ \Gamma $-covering $G^{ \alpha } $ of $G$ is 
connected. 
Then the Kirchhoff index of $G^{ \alpha}$ is
\[
Kf^z (G^{ \alpha } )=| \Gamma | Kf^z (G)-2(m-n)n(| \Gamma |^2 -| \Gamma |) 
\]
\[
+2| \Gamma |(m-n) \sum^k_{i=2} f_i 
\frac{d}{dt} \log \det ( {\bf I}_{f_i n} -t \sum_{g \in \Gamma } { \rho }_i (g) \bigotimes {\bf A}_{g} 
+ t^2 ( f_i \circ {\bf Q} )) |_{t=1} . 
\]
\end{corollary}

{\bf Acknowledgments}

We would like to thank the referees for many useful suggestions and comments. 
The first author is partially supported by the Grant-in-Aid for Scientific Research (C) of Japan
Society for the Promotion of Science (Grant No. 16K05187). 
The second author is partially supported by the Grant-in-Aid for Scientific Research (C) of Japan
Society for the Promotion of Science (Grant No. 16K05249). 
The third author is  partially supported by the Grant-in-Aid for Young Scientists (B) of Japan Society 
for the Promotion of Science (Grant No. 26400001). 
The forth author is partially supported by the Grant-in-Aid for Scientific Research (C) of Japan
Society for the Promotion of Science (Grant No. 15K04985).

\end{document}